\documentclass{lms}
\usepackage{amssymb}
\usepackage{lmsobit}
\usepackage{graphicx}
\usepackage{rotating}

\newnumbered{notation}[theorem]{Notation}

\newnumbered{assertion}{Assertion}    
\newnumbered{conjecture}{Conjecture}  
\newnumbered{definition}{Definition}
\newnumbered{hypothesis}{Hypothesis}
\newnumbered{remark}{Remark}
\newnumbered{note}{Note}
\newnumbered{observation}{Observation}
\newnumbered{problem}{Problem}
\newnumbered{question}{Question}
\newnumbered{algorithm}{Algorithm}
\newnumbered{example}{Example}

\title[Obituary]{OBITUARY}
\obitfor[Graham Everest]{{\large Graham Everest, 1957--2010}}

\begin{document}
\maketitle

\vspace{-6mm}

\begin{figure}[!h]
\begin{center}
%
\includegraphics[height=176pt]{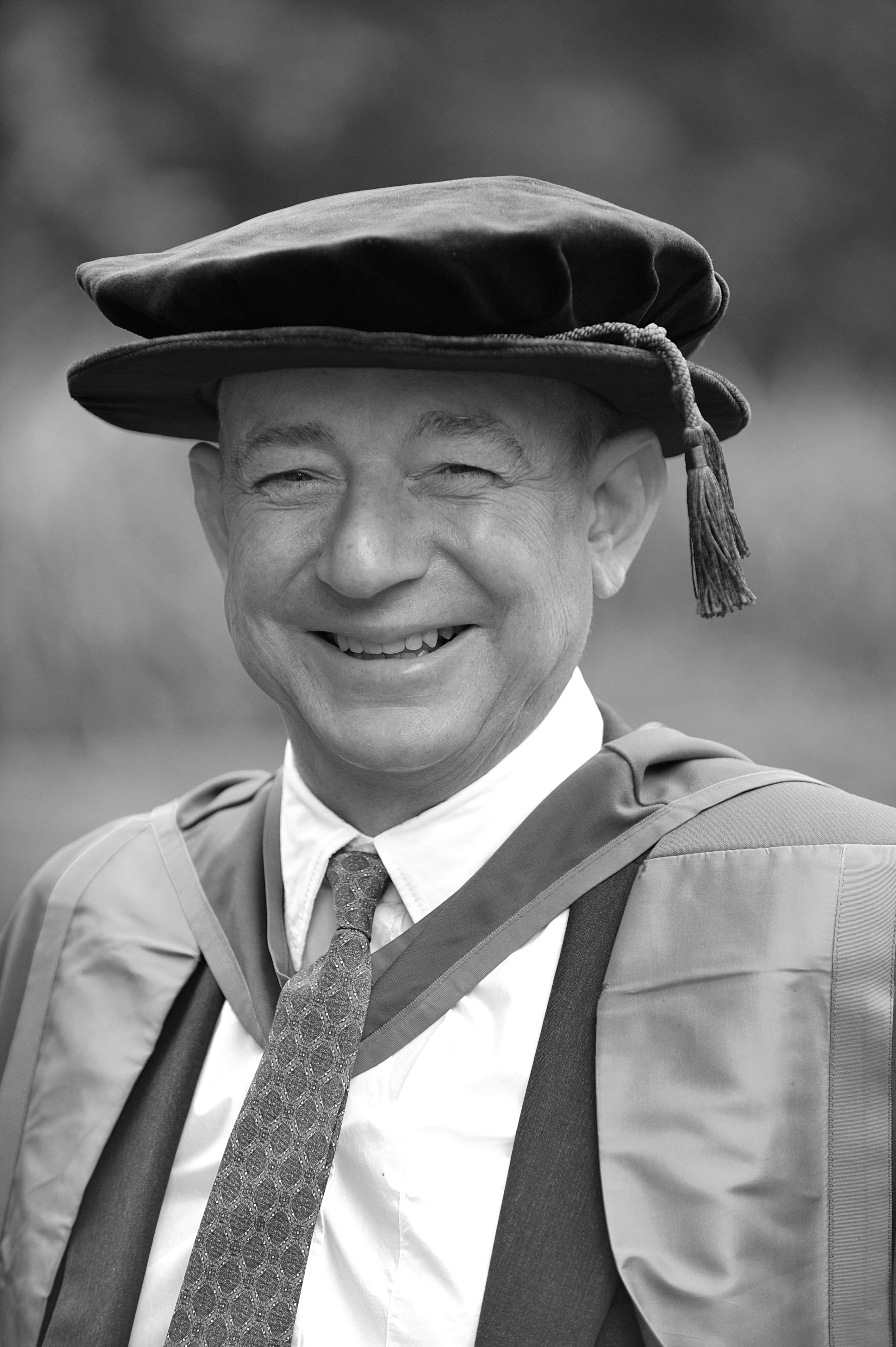}%
\begin{turn}{90}{\footnotesize\rule{0pt}{3mm}\tiny\LMScopyright University of East Anglia}\end{turn}%
\hspace*{-4mm}
\end{center}
\end{figure}


\begin{quote}
``He heals the brokenhearted and binds up their wounds.
He determines the number of the stars and calls them each by name''
Psalm 147, v. 3--4 (NIV).
\end{quote}


\section{Introduction}

\nocite{*}

\noindent Graham Everest was an inspirational mathematician, who touched the mathematical
lives of a great many more than the thirty with whom he published.
He collaborated widely, and those who worked with him enjoyed the
joy he took in learning and doing mathematics. His creative leaps could generate
mistakes, but these were often of the fruitful variety, triggering new ideas and
ways forward. While his primary research area was in algebraic number theory,
he was always interested in interactions between different parts of mathematics,
and his research publications touched on aspects of analytic and computational number
theory, dynamical systems, and logic. He was an enthusiastic user of computational
methods to inform theoretical thinking and to enhance his approach to teaching.
Alongside conventional research publications, the exposition of mathematics
mattered greatly to him. While it sadly came after his death, the fact that one of his
publications~\citep{MR2826450} won the 2012 Lester Ford Award from the Mathematical
Association of America would have
brought him huge pleasure. He was a dedicated and enthusiastic teacher, and deservedly
won a UEA Excellence in Teaching Award from student
nominations in 2005.

\section{Personal life}

As a person Graham was warm, generous, and possessed of a large hinterland.
This ranged from playing bridge, amateur dramatics, interests in world music,
the writings of Carl Jung, religious mysticism and cake-baking, to arranging wine
and poetry evenings. His Christianity eventually saw him move away from
Evangelical free churches and an --- at times controversial --- advocacy of
young earth creationism~\citep{7days}, towards the Church of England and a humble
and questioning approach to all such matters.
Eventually Graham's faith led to his
undergoing theological training,
initially as a lay reader, and later to his being ordained
as a non-Stipendiary Minister
in~2006
at a ceremony in Norwich Cathedral. Parishioners from Colney and Cringleford
in Norwich appreciated the warmth, wit, and love he brought to his
vocation, in which he felt a particular calling to the ministry of
encouragement.

The year~2008 brought enormous difficulties to Graham's life.
At the start of the year he received the devastating
diagnosis of an aggressive and inoperable form of cancer.
Later his own at times complicated personal struggles led to him
leaving his role in the Church.
He faced death much as he had faced life, with honesty,
candour, wit, and passion. Shortly before he passed away
at home, he
was reconciled with the Church and reinstated as a priest.

Graham remained mathematically active throughout his illness,
and continued to teach while he was physically able to. One of
his last publications~\citep{MR2826450} was dedicated to
`Mulbarton Ward and the Weybourne Day Unit in Norwich',
where he had received medical and hospice care respectively.

Graham is survived by his children James, Philip, and Rebekah, and his wife Sue.

\section{Mathematical career}

Graham was born into a working-class family where University was not
seen as an automatic path to follow,
in Southwick, West Sussex. His early childhood was happy,
and his playful personality developed through inventive games played with his
older sister Jan. Primary school was not an easy time for Graham, perhaps
as a result of boredom. His mathematical talent was spotted at High School,
and Graham was the first person at his school to take an S-level examination
in Mathematics. More interested in playing bridge than in
his formal studies at school,
he dropped an A-level in order to devote more time to this.
He followed his father's advice, and started an auditing apprenticeship --
but this failed to fulfil him intellectually. Happily,
Graham's Mathematics teacher
helped him apply to study Mathematics in the next academic year.

Bedford College (now part of Royal Holloway, University of London),
where Graham started his degree in~1977,
expanded Graham's horizons in many ways. He embraced literature, art,
and music, starting several life-long interests. One of several
great friends he made, who later went on to become a monk, shared his Christian
faith with Graham, and the outworking of this faith was to influence
all aspects of Graham's later life, and to have a profound influence on
members of his family and many of his friends.
Graham's mathematical powers were clear, and he won both
first and second year
prizes, later
explaining to his wife Sue that the missing
third year prize
had not been won because
it did not exist.
He went on to study for a doctorate
at King's College London
under the supervision of Colin Bushnell,
and during this period was both influenced and supported by
Martin Taylor and Albrecht Fr{\"o}hlich. In the Summer of~1983 he
completed his PhD thesis~\citep{thesis}
and married Sue, who had by then also graduated from
Bedford College. Graham turned down
a temporary lectureship at Sheffield (who were
kind enough to release him)
in order to take up a permanent lectureship at UEA.

During his long career at UEA, Graham was always
energetically dedicated to both education and research.
For some years in the~1980s the flame of pure mathematics research
at UEA was largely kept alive by the efforts of Graham and
Alan Camina. Graham took great
pleasure in the steady growth in the research power of the UEA
Mathematics Department that
took place thereafter, and was openly delighted by
the vibrant atmosphere that emerged.

\section{Mathematical research}

Graham's research outputs appeared in the form of some~70
research
and research-expository papers and three monographs~\citep{MR1700272},~\citep{MR1990179},
and~\citep{MR2135478}. While always remaining number-theoretical,
his work may be roughly categorized into three main areas:
Diophantine analysis, dynamical systems, and recurrence sequences.

\subsection{Diophantine analysis and counting problems}\label{sectionDiophantine}

Graham's early interests lay in using analytic and Diophantine
methods to understand counting and distribution problems.
In his first papers he studied a conjecture of Bushnell~\citer{MR0546665},
using Schmidt's subspace theorem to study
the distribution of normal integral generators in
number fields~\citep{MR0695353},~\citep{MR0713379}.
He saw how one could employ Alan Baker's extremely powerful work on
Diophantine approximation~\citer{MR0220680}
in the investigation of quite complex algebraic constraints. This particular
combination was unprecedented at the time, and already exhibited the originality and insight
of Graham's later work.
He went on to use analytic methods in the same area~\citep{MR0771814},
setting the tone for much of his later work: bringing deep
Diophantine results and an analytic toolbox to bear on
counting or growth problems.

Motivated by the work of Evertse~\citer{MR0766298} on
the~$S$-unit theorem, Graham used analytic methods
of the Hardy--Littlewood type
to study the distribution
of~$N(\mathbf{x})=\prod_v\vert c_0x_0+\cdots+c_nx_n\vert_v$ for~$x_i$
running through~$S$-units in a number field~\citep{MR0996323}.
He also studied norm-form equations from a similar
perspective~\citep{MR0957247},~\citep{MR1165359},~\citep{MR1438600}.

\subsection{Dynamical systems of algebraic origin}

I first met Graham across an interview table in 1991. Despite
the constraints of the setting, it was quickly apparent that
we shared an interest in the remarkable, then recent, work
of Schlickewei and others on the subspace theorem,
and the resulting insights into the solution of~$S$-unit
equations
in fields of characteristic zero~\citer{MR0948309},~\citer{MR1119694},~\citer{MR1069241}.
In algebraic dynamics
this had led directly to the result that a mixing
action of~$\mathbb Z^d$ by automorphisms of a compact
connected group is mixing of all orders~\citer{MR1193598},
and for Graham this unexpected instance of number-theoretical
results finding applications in dynamical systems was fascinating.
This began a conversation between us that lasted almost twenty
years, each bringing in tools and techniques from our own
areas of mathematics.

Algebraic dynamics is also a setting in which the
Mahler measure
\[
m(f)=\int_0^1\cdots\int_0^1\log\vert f({\rm e}^{2\pi{\rm i}s_1},\dots,{\rm e}^{2\pi{\rm i}s_d})\vert
{\rm d}s_1\cdots{\rm d}s_d
\]
of a polynomial~$f\in\mathbb{Z}[u_1^{\pm1},\dots,u_d^{\pm1}]$ arises
naturally as the topological entropy of a~$d$-dimensional dynamical
system~\citer{MR1062797} associated to the polynomial~$f$. In the case~$d=1$ the most natural
proof exploits properties of the adele ring~$\mathbb Q_{\mathbb A}$,
and in the simplest case of linear polynomials the notion of Diophantine
height appears in a particularly transparent way~\citer{MR0961739}.
If~$d=1$ and~$f(u_1)=bu_1-a$ then the associated dynamical system
is given by the map~$\alpha$ dual to~$x\mapsto\frac{a}{b}x$ on the character group of~$\mathbb Z[\frac{1}{b}]$,
and the topological entropy is given by
\[
h(\alpha)=\sum_{p\le\infty}\log\max\{0,\vert\textstyle\frac{a}{b}\vert_p\}
=\log\max\{\vert a\vert,\vert b\vert\}=m(f).
\]
Once again Graham was
intrigued to find something he had studied in a number-theoretical
setting making an appearance in dynamical systems.
Two aspects of Mahler measure
became of particular interest to Graham --- its special values, which later
came to have interpretations in terms of periods
of mixed motives~\citer{MR1415320}, and the classical
problem of Lehmer~\citer{MR1503118} asking if
\[
\inf\{m(f)\mid m(f)>0\}>0.
\]
This rich circle of questions, a mixture of Diophantine
analysis and arithmetic geometry, triggered a long interest in the
interplay between dynamical systems and arithmetic,
leading to several research projects and the monograph~\citep{MR1700272}
on the relationship between heights of algebraic numbers and
topological entropy in algebraic dynamics.

One theme that Graham pursued in this area concerned
growth in periodic orbits for automorphisms of solenoids,
which relates directly to questions about the
arithmetic of terms of recurrence sequences and linked
to his long-standing interest in the arithmetic of
linear and bi-linear recurrence sequences.
These questions are of interest in algebraic
dynamics because they arise in any attempt to
understand typical or generic
compact group automorphisms.
The
simplest examples boil down to questions of the
following shape. What can be said about growth
in expressions like
\[
(2^n-1)\prod_{p\in S}\vert 2^n-1\vert_p
\]
as~$n\to\infty$ for various subsets~$S$ of
the set of all primes~$\mathbb{P}$?
Trying to develop a good understanding of what can
be said here, and how much averaging is needed to smooth
out the Mersenne-like quirks in the prime factorization
of linear recurrence sequences, held Graham's attention
on and off for several years, and eventually
led to a rather complete understanding of the
two extreme cases~$\vert S\vert<\infty$ in~\citep{MR2339472}
and~$\vert\mathbb{P}\setminus S\vert<\infty$ in~\citep{MR2550149}.
The latter case brought Graham particular pleasure as it
involved Dirichlet series and analytic issues relating
to an asymptotic counting problem in a novel setting,
often with such poor analytic behaviour on the critical line
that Tauberian theorems could not be applied.
The large middle ground ---
the uncountable collection of infinite sets
of primes with infinite complement --- saw
some early partial results~\citep{MR1461206}
leading eventually to constructions of families
of examples exhibiting continua of different
orbit growth
rates for a suitably averaged measure of
orbit growth~\citer{baierjaideestevensward}.

\subsection{Recurrence sequences, heights, and logic}

No mathematician is unaware of open problems concerned with
the appearance of
primes in the sequence~$(2^n-1)$ (`Mersenne' primes); only slightly less well-known is the
well-understood appearance of primitive divisors in the same
sequence (a primitive divisor of a term is a prime divisor that
does not divide any earlier term).
Many mathematicians might find the briar patch
of notorious difficulties
surrounding Mersenne primes, or the over-manicured
and well-trodden garden surrounding the question of primitive
divisors in sequences like~$(a^n-b^n)$, with its
frequently re-proved and definitive Zsigmondy--Bang theorem~\citer{MR1546236}, less than appealing.
It was typical of Graham that in both cases he saw
rich territory for new exploration and extension,
and he was particularly enthused by the way in which
Bilu, Hanrot and Voutier~\citer{MR1863855} were able to use
modern Diophantine results to
solve a century-old problem, proving that
the~$n$th term of any Lucas or Lehmer sequence has a primitive divisor
for~$n>30$. This paper, with its sophisticated use of deep
Diophantine results and technical skill to produce a result
that was startling both in its uniformity and in its
highly effective bounds, was a great inspiration to Graham.

Graham studied the arithmetic of linear recurrence sequences
from many different perspectives, and some of
these questions were pursued by doctoral students under
his supervision. The Zsigmondy--Bang theorem involves
some key estimates requiring an understanding of the
rate of growth in a linear recurrence sequence,
and the delicate issues arising here may already be
seen in Lehmer--Pierce sequences~$(\Delta_n(f))$, where~$f(x)=\prod_{i=1}^d(x-\alpha_i)$
is a monic integer polynomial with no cyclotomic factor,
and~$\Delta_n(f)=\prod_{i=1}^{d}\vert\alpha_i^n-1\vert$.
If no~$\alpha_i$ has unit modulus, then the
growth rate is clear,
while if~$f$ has zeros of
unit modulus then a {\it deus ex machina} like
Baker's theorem is required to know that
\[
\lim_{n\to\infty}\frac{1}{n}\log\Delta_n=\sum_{\vert\alpha_i\vert>1}\log\vert\alpha_i\vert
=m(f),
\]
the logarithmic Mahler measure or height of~$f$.
The question of prime appearance remains inaccessible for
linear recurrence sequences, triggering one of Graham's excursions
into computational work. Writing~$(n_j)$ for the sequence of indices
for which~$\Delta_{n_j}$ is prime, he used Baker's theorem
to show that~$(n_j)$ has only finitely many composite terms~\citep{MR1783409}, and
built on Wagstaff's heuristics to argue that~$j/\log\log\Delta_{n_j}$ should
converge as~$j\to\infty$ to a quantity controlled by~$m(f)$,
adding further numerical evidence to the original insight of Lehmer~\citer{MR1503118}.
Some of these ideas found more rigorous
outlet in work of his students, notably that of Flatters~\citer{MR2468478}, where
uniform bounds for the appearance of primitive divisors
in sequences associated to real quadratic units are found.

Graham was also interested in elliptic or bi-linear recurrence
sequences.
For an elliptic curve~$E$ defined over~$\mathbb Q$ and
given in Weierstrass form, and a non-torsion point~$P\in E(\mathbb{Q})$,
there is an associated integer sequence~$(B_n)$ defined by~$[n]P=(\frac{A_n}{B_n^2},\frac{C_n}{B_n^3})$.
Here the work of Silverman~\citer{MR0961918} showing that all but
finitely many terms of an elliptic divisibility sequence have
a primitive divisor proved to be inspirational.
Graham formulated a conjectural view of
the arithmetic properties of these sequences, including the
striking suggestion that the number of prime terms in~$(B_n/B_1)$
should be bounded uniformly. He was able to establish
many special cases, including strong results
conditional on Lang's height conjecture~\citep{MR2429645},~\citep{MR2045409};
generalizations to Somos sequences associated to
sequences like~$([n]P+Q)$ in~\citep{MR2200959};
generalizations to Siegel and Hall theorems~\citep{MR2548983};
and results on
primitive divisors~\citep{MR2486632}, including highly effective bounds across
certain families of curves~\citep{MR2220263}.

In both the linear and the elliptic theory, a critical role
is played by notions of height. In the linear case, this
may be expressed as a Mahler measure, and in the elliptic
case as a Neron--Tate height. Pursuing a broad thematic analogy
between the two theories, Graham introduced in~\citep{MR1373343} an elliptic
Mahler measure~$m_E(F)$ as a sum of
local integrals,
using an interesting integral representation for the canonical local height on~$E$.
In particular, this gives rise to a beautiful elliptic analogue of Jensen's theorem, and
the result that~$m_E(F)=0$ if and only if the roots of~$F$ are the~$x$-coordinates
of division points (a direct elliptic analogue of Kronecker's lemma in the classical case).

The possibility of relating Hilbert's 10th problem for~$\mathbb Q$ to
arithmetic properties of elliptic divisibility sequences
(see Poonen~\citer{poonen}) fascinated Graham,
leading to an elegant argument~\citep{MR2480276}
building on Poonen's work
to find a partition~$S\sqcup T$ of the primes into
recursive sets with the
property that Hilbert's 10th problem is undecidable
for both~$\mathbb Z_S$ and~$\mathbb Z_T$
(later generalized to number fields in~\citep{MR2915472}).
How Graham came to be interested in these questions
highlights several aspects of the way he worked.
Gunther Cornelissen, following a talk by Thanases Pheidas in Gent in~1999
raising the possibility of attacking the problem using arithmetic
properties of elliptic divisibility sequences, started exploring some
of the literature on these sequences, leading to the
results in~\citer{MR2212121}.
For quite different reasons, Graham had been developing
a web site collecting material and his own thoughts
on elliptic divisibility sequences, bringing together
quite disparate strands of thoughts growing from
Morgan Ward's classical treatment~\citer{MR0023275}.
Gunther came across this web site while
researching the topic.
This led to Graham attending a
mini-workshop in Oberwolfach on Hilbert's Tenth Problem, and
to collaboration with some of the participants, including
Alexandra Shlapentokh~\citep{MR2915472} and Kirsten Eisentr{\"a}ger~\citep{MR2480276},~\citep{MR2915472}.
Using the web to find another mathematical environment with new concepts
to learn, people to meet, and friends to make, where his knowledge of elliptic divisibility sequences
could be applied, brought him great pleasure in the last few years of
his life.

The impact Graham had as a person on several different mathematical communities
was reflected in two international conferences.
The first --- `Diverse faces of arithmetic' --- which he much enjoyed, was a retirement conference at UEA in December 2009, bringing
together people from number theory, dynamical systems, integrable systems, and logic.
The second, `Definability in Number Theory' at the University of Gent in September 2010,
was dedicated to his memory and
brought together
mathematicians
interested in the
question of which sets and structures can be defined or interpreted in the existential or first-order theory of rings and fields.

\section{Research students}

Graham was a dedicated and thoughtful supervisor. The breadth
of his interests and his determination to continue to
supervise graduate students as his illness progressed
meant
that joint supervisions arose
naturally, and both Shaun Stevens and I enjoyed several joint supervisory
experiences with him.

\medskip
\centerline{\sc Research students supervised or jointly supervised by Graham Everest}
\smallskip
\begin{enumerate}
\item[] Alice Miller, PhD  1988, {\it`Effective subspace theorems for function fields'}
\item[] Br{\'{\i}}d N{\'{\i}} Fhlath{\'u}in, PhD 1995, {\it`Mahler's measure on Abelian varieties'}
\item[] Peter Panayi, PhD 1995, {\it`Computation of Leopoldt's $p$-adic regulator'}
\item[]	Vijay Chothi, PhD 1996 (with Tom Ward), {\it`Periodic points in~$S$-integer dynamical systems'}
\item[] Paola D'Ambros, PhD 2000 {\it`Algebraic dynamics in positive characteristic'}
\item[] Christian R{\"o}ttger, PhD 2000, {\it`Counting problems in algebraic number theory'}
\item[] Peter Rogers, MPhil 2003, {\it`Computational aspects of elliptic curves'}
\item[] Patrick Moss, PhD 2003, {\it`The arithmetic of realizable sequences'}
\item[] Victoria Stangoe, PhD 2004 (with Tom Ward), {\it`Orbit counting far from hyperbolicity'}
\item[] Helen King, PhD 2005 (with Tom Ward), {\it`Prime appearance in elliptic divisibility sequences'}
\item[] Jonathan Reynolds, PhD 2008 (with Shaun Stevens), {\it`Extending Siegel's theorem for elliptic curves'}
\item[] Ouamporn Phuksuwan, PhD 2009 (with Shaun Stevens), {\it`The uniform primality conjecture for the twisted Fermat cubic'}
\item[] Anthony Flatters, PhD 2010 (with Tom Ward), {\it`Arithmetic properties of recurrence sequences'}
\end{enumerate}

\begin{acknowledgements}
I am grateful to Sue and James Everest for providing
information about Graham's early years, to Colin Bushnell for
comments on Section~\ref{sectionDiophantine}, to Gunther
Cornelissen for comments on Hilbert's Tenth Problem, and to
Maresa Padmore for locating the photograph of Graham being made
an Emeritus Professor at UEA. I am also grateful for comments
from Alan Camina, Kirsten Eisentr{\"a}ger, Anthony Flatters,
Thanases Pheidas, Christian R{\"o}ttger, and David Stevens.
Particular thanks are due to Shaun Stevens with whom I had
several discussions about an obituary, and to Anish Ghosh,
Shaun Stevens, and Sanju Velani for organising Graham's
retirement conference which meant so much to him.
\end{acknowledgements}

\begin{anglerefs}{99}{References}

\bibitem{baierjaideestevensward} S.~Baier, S.~Jaidee,
    S.~Stevens and T.~Ward, `Automorphisms with exotic orbit
    growth', \emph{Acta Arithmetica}, \textbf{158}, no.~2 (2013), 173--197.

\bibitem{MR0220680}
A.~Baker,
`Linear forms in the logarithms of algebraic numbers. I, II, III, IV',
\emph{Mathematika} \textbf{13} (1966), 204--216; {\it ibid.} \textbf{14} (1967), 102--107; {\it ibid.} \textbf{14} (1967), 220--228;
{\it ibid.} \textbf{15} (1968), 204--216.

\bibitem{MR1863855}
Yu.~Bilu, G.~Hanrot and P.~M.~Voutier,
`Existence of primitive divisors of Lucas and Lehmer numbers,
(with an appendix by M. Mignotte)',
\emph{J. Reine Angew. Math.} \textbf{539} (2001), 75--122.

\bibitem{MR0546665}
C.~J.~Bushnell,
`Norm distribution in Galois orbits',
\emph{J. Reine Angew. Math.} \textbf{310} (1979), 81--99.

\bibitem{MR2212121}
G.~Cornelissen, T.~Pheidas and K.~Zahidi,
`Division-ample sets and the Diophantine problem for rings of integers',
\emph{J. Th{\'e}or. Nombres Bordeaux} \textbf{17} (2005), no.~3, 727--735.

\bibitem{MR1415320}
C.~Deninger,
`Deligne periods of mixed motives, $K$-theory and the entropy of certain $\mathbb Z^n$-actions',
\emph{J. Amer. Math. Soc.} \textbf{10} (1997), no.~2, 259--281.

\bibitem{MR0766298}
J.-H.~Evertse,
`On sums of $S$-units and linear recurrences',
\emph{Compositio Math.} \textbf{53} (1984), no.~2, 225--244.

\bibitem{MR0948309}
J.-H.~Evertse and K.~Gyory,
`On the numbers of solutions of weighted unit equations',
\emph{Compositio Math.} \textbf{66} (1988), no.~3, 329--354.

\bibitem{MR2468478}
A.~Flatters,
`Primitive divisors of some Lehmer--Pierce sequences',
\emph{J. Number Theory} \textbf{129} (2009), no.~1, 209--219.

\bibitem{MR1503118}
D.~H.~Lehmer,
`Factorization of certain cyclotomic functions',
\emph{Ann. of Math. (2)} \textbf{34} (1933), no.~3, 461--479.

\bibitem{MR1062797}
D.~A.~Lind, K.~Schmidt and T.~Ward,
`Mahler measure and entropy for commuting automorphisms of compact groups',
\emph{Invent. Math.} \textbf{101} (1990), no.~3, 593--629.

\bibitem{MR0961739}
D.~A.~Lind and T.~Ward,
`Automorphisms of solenoids and $p$-adic entropy',
\emph{Ergodic Theory Dynam. Systems} \textbf{8} (1988), no.~3, 411--419.

\bibitem{poonen}
B.~Poonen,
`Hilbert's Tenth Problem and Mazur's Conjecture
for large subrngs of $\mathbb Q$',
\emph{J. Amer. Math. Soc.}, \textbf{16} (2002), no.~4, 981--990.

\bibitem{MR1119694}
A.~J.~van der Poorten and H.~P.~Schlickewei, H. P.,
`Additive relations in fields',
\emph{J. Austral. Math. Soc. Ser.~A} \textbf{51} (1991), no.~1, 154--170.

\bibitem{MR1069241}
H.~P.~Schlickewei,
`$S$-unit equations over number fields',
\emph{Invent. Math.} \textbf{102} (1990), no.~1, 95--107.

\bibitem{MR1193598}
K.~Schmidt and T.~Ward,
`Mixing automorphisms of compact groups and a theorem of Schlickewei',
\emph{Invent. Math.} \textbf{111} (1993), no.~1, 69--76.

\bibitem{MR0961918}
J.~H.~Silverman,
`Wieferich's criterion and the $abc$-conjecture',
\emph{J. Number Theory} \textbf{30} (1988), no.~2, 226--237.

\bibitem{MR0023275}
M.~Ward,
`Memoir on elliptic divisibility sequences',
\emph{Amer. J. Math.} \textbf{70} (1948), 31--74.

\bibitem{MR1546236}
K.~Zsigmondy,
`Zur Theorie der Potenzreste',
\emph{Monatsh. Math. Phys.} \textbf{3} (1892), no.~1, 265--284.

\end{anglerefs}

\begin{squarerefs}{999}{Publications of Graham Everest}

\bibitem{MR0695353}
`The distribution of normal integral generators', in
  \emph{Seminar on {N}umber {T}heory, 1981/1982}, pp.~Exp. No. 43, 8 (Univ.
  Bordeaux I, Talence, 1982).

\bibitem{thesis}
`The distribution of normal integral generators in tame extensions of~$\mathbb{Q}$',
PhD thesis, King's College London (1983).

\bibitem{MR0713379}
`Diophantine approximation and the distribution of normal
  integral generators', \emph{J. London Math. Soc. (2)} \textbf{28} (1983),
  no.~2, 227--237.

\bibitem{MR0816482}
`Independence in the distribution of normal integral
  generators', \emph{Quart. J. Math. Oxford Ser. (2)} \textbf{36} (1985),
  no.~144, 405--412.

\bibitem{MR0771814}
`Diophantine approximation and {D}irichlet series', \emph{Math.
  Proc. Cambridge Philos. Soc.} \textbf{97} (1985), no.~2, 195--210.

\bibitem{MR829584}
(with {\sc A.~R. Camina} and {\sc T.~M. Gagen})  `Enumerating nonsoluble
  groups---a conjecture of {J}ohn {G}. {T}hompson', \emph{Bull. London Math.
  Soc.} \textbf{18} (1986), no.~3, 265--268.

\bibitem{MR0824441}
`The divisibility of normal integral generators', \emph{Math.
  Z.} \textbf{191} (1986), no.~3, 397--404.

\bibitem{MR0878134}
`Galois generators and the subspace theorem', \emph{Manuscripta
  Math.} \textbf{57} (1987), no.~4, 451--467.

\bibitem{MR0882289}
`Angular distribution of units in abelian group rings---an
  application to {G}alois-module theory', \emph{J. Reine Angew. Math.}
  \textbf{375/376} (1987), 24--41.

\bibitem{MR0993111}
`Some meromorphic functions associated to the {$S$}-unit
  equation', in \emph{S\'eminaire de {T}h\'eorie des {N}ombres, 1987--1988
  ({T}alence, 1987--1988)}, pp.~Exp.\ No.\ 12, 10 (Univ. Bordeaux I, Talence,
  1988).

\bibitem{MR0957247}
`A ``{H}ardy-{L}ittlewood'' approach to the norm form
  equation', \emph{Math. Proc. Cambridge Philos. Soc.} \textbf{104} (1988),
  no.~3, 421--427.

\bibitem{MR1007893}
`Units in abelian group rings and meromorphic functions',
  \emph{Illinois J. Math.} \textbf{33} (1989), no.~4, 542--553.

\bibitem{MR0996323}
`A ``{H}ardy-{L}ittlewood'' approach to the {$S$}-unit
  equation', \emph{Compositio Math.} \textbf{70} (1989), no.~2, 101--118.

\bibitem{MR0987018}
`Root numbers---the tame case', in \emph{Representation theory
  and number theory in connection with the local {L}anglands conjecture
  ({A}ugsburg, 1985)}, in \emph{Contemp. Math.} \textbf{86}, pp.~109--116
  (Amer. Math. Soc., Providence, RI, 1989).

\bibitem{MR1072046}
`A new invariant for tame, abelian extensions', \emph{J. London
  Math. Soc. (2)} \textbf{41} (1990), no.~3, 393--407.

\bibitem{MR1062335}
`Counting the values taken by sums of {$S$}-units', \emph{J.
Number Theory} \textbf{35} (1990), no.~3, 269--286.

\bibitem{MR1058239}
`The {$S$}-unit equation and {D}irichlet series', in
  \emph{Number theory, {V}ol.\ {II} ({B}udapest, 1987)}, in \emph{Colloq. Math.
  Soc. J\'anos Bolyai} \textbf{51}, pp.~659--669 (North-Holland, Amsterdam,
  1990).

\bibitem{MR1185339}
(with {\sc T.~M. Gagen})  `Measures associated to the inverse regulator
  of a number field', \emph{Arch. Math. (Basel)} \textbf{59} (1992), no.~5,
  420--426.

\bibitem{MR1165549}
`Applications of the {$p$}-adic subspace theorem', in
  \emph{{$p$}-adic methods and their applications}, in \emph{Oxford Sci.
  Publ.}, pp.~33--56 (Oxford Univ. Press, New York, 1992).

\bibitem{MR1175691}
`Addendum: ``{O}n the solution of the norm-form equation''',
  \emph{Amer. J. Math.} \textbf{114} (1992), no.~4, 787--788.

\bibitem{MR1184758}
`On the canonical height for the algebraic unit group',
  \emph{J. Reine Angew. Math.} \textbf{432} (1992), 57--68.

\bibitem{MR1165359}
`On the solution of the norm-form equation', \emph{Amer. J.
  Math.} \textbf{114} (1992), no.~3, 667--682.

\bibitem{MR1179007}
`{$p$}-primary parts of unit traces and the {$p$}-adic
  regulator', \emph{Acta Arith.} \textbf{62} (1992), no.~1, 11--23.

\bibitem{MR1164774}
`Uniform distribution and lattice point counting', \emph{J.
  Austral. Math. Soc. Ser. A} \textbf{53} (1992), no.~1, 39--50.

\bibitem{MR1225955}
(with {\sc J.~H. Loxton}) `Counting algebraic units with bounded
  height', \emph{J. Number Theory} \textbf{44} (1993), no.~2, 222--227.

\bibitem{MR1269287}
`An asymptotic formula implied by the {L}eopoldt conjecture',
  \emph{Quart. J. Math. Oxford Ser. (2)} \textbf{45} (1994), no.~177, 19--28.

\bibitem{MR1250998}
`Corrigenda to: ``{U}niform distribution and lattice point
  counting'' [{J}.\ {A}ustral.\ {M}ath.\ {S}oc.\ {S}er.\ {A} {\bf 53} (1992),
  no.\ 1, 39--50; {MR}1164774 (93i:11114)]', \emph{J. Austral. Math. Soc. Ser.
  A} \textbf{56} (1994), no.~1, 144.

\bibitem{MR1316819}
`On the proximity of algebraic units to divisors', \emph{J.
  Number Theory} \textbf{50} (1995), no.~2, 233--250.

\bibitem{MR1335283}
`On the {$p$}-adic integral of an exponential polynomial',
\emph{Bull. London Math. Soc.} \textbf{27} (1995), no.~4, 334--340.

\bibitem{MR1317513}
`Mean values of algebraic linear forms', \emph{Proc. London
  Math. Soc. (3)} \textbf{70} (1995), no.~3, 529--555.

\bibitem{MR1332880}
`The mean value of a sum of {$S$}-units', \emph{J. London Math.
  Soc. (2)} \textbf{51} (1995), no.~3, 417--428.

\bibitem{MR1313121}
`Estimating {M}ahler's measure', \emph{Bull. Austral. Math.
  Soc.} \textbf{51} (1995), no.~1, 145--151.

\bibitem{MR1380647}
(with {\sc V.~Chothi} and {\sc T.~Ward})  `Oriented local entropies for expansive
  actions by commuting automorphisms', \emph{Israel J. Math.} \textbf{93}
  (1996), 281--301.

\bibitem{MR1382488}
(with {\sc I.~E. Shparlinski}) `Divisor sums of generalised exponential
  polynomials', \emph{Canad. Math. Bull.} \textbf{39} (1996), no.~1, 35--46.

\bibitem{MR1373343}
(with {\sc B.~N. Fhlath{\'u}in})  `The elliptic {M}ahler measure',
  \emph{Math. Proc. Cambridge Philos. Soc.} \textbf{120} (1996), no.~1, 13--25.

\bibitem{MR1461206}
(with {\sc V.~Chothi} and {\sc T.~Ward})  `{$S$}-integer dynamical systems: periodic
  points', \emph{J. Reine Angew. Math.} \textbf{489} (1997), 99--132.

\bibitem{MR1401740}
(with {\sc A.~J. van~der Poorten})  `Factorisation in the ring of
  exponential polynomials', \emph{Proc. Amer. Math. Soc.} \textbf{125} (1997),
  no.~5, 1293--1298.

\bibitem{MR1438600}
(with {\sc K.~Gy{\H{o}}ry})  `Counting solutions of decomposable form
  equations', \emph{Acta Arith.} \textbf{79} (1997), no.~2, 173--191.

\bibitem{MR1678099}
(with {\sc T.~Ward})  `A dynamical interpretation of the global canonical
  height on an elliptic curve', \emph{Experiment. Math.} \textbf{7} (1998),
  no.~4, 305--316.

\bibitem{MR1666050}
(with {\sc C.~Pinner}) `Bounding the elliptic {M}ahler measure. {II}',
  \emph{J. London Math. Soc. (2)} \textbf{58} (1998), no.~1, 1--8.

\bibitem{MR1646573}
`Measuring the height of a polynomial', \emph{Math. Intelligencer}
  \textbf{20} (1998), no.~3, 9--16.

\bibitem{MR1646051}
`Counting generators of normal integral bases', \emph{Amer. J.
  Math.} \textbf{120} (1998), no.~5, 1007--1018.

\bibitem{MR1700272}
(with {\sc T.~Ward}) \emph{Heights of polynomials and entropy in algebraic
  dynamics}, in \emph{Universitext} (Springer-Verlag London Ltd., London,
  1999).

\bibitem{MR1485471}
(with {\sc I.~E. Shparlinski}) `Counting the values taken by algebraic
  exponential polynomials', \emph{Proc. Amer. Math. Soc.} \textbf{127} (1999),
  no.~3, 665--675.

\bibitem{MR1688486}
`On the elliptic analogue of {J}ensen's formula', \emph{J.
  London Math. Soc. (2)} \textbf{59} (1999), no.~1, 21--36.

\bibitem{MR1703208}
`Explicit local heights', \emph{New York J. Math.} \textbf{5}
  (1999), 115--120.

\bibitem{MR1778843}
(with {\sc P.~D'Ambros, R.~Miles}, and {\sc T.~Ward})  `Dynamical systems arising
  from elliptic curves', \emph{Colloq. Math.} \textbf{84/85} (2000), no.~1, 95--107.
\newblock Dedicated to the memory of Anzelm Iwanik.

\bibitem{MR1783650}
(with {\sc C.~Pinner})  `Corrigendum: ``{B}ounding the elliptic {M}ahler
  measure. {II}'' [{J}. {L}ondon {M}ath. {S}oc. (2) {\bf 58} (1998), no. 1,
  1--8; {MR}1666050 (2000a:11099)]', \emph{J. London Math. Soc. (2)}
  \textbf{62} (2000), no.~2, 640.

\bibitem{MR1783409}
(with {\sc M.~Einsiedler} and {\sc T.~Ward})  `Primes in sequences associated to
  polynomials (after {L}ehmer)', \emph{LMS J. Comput. Math.} \textbf{3} (2000),
  125--139.

\bibitem{MR1800354}
(with {\sc T.~Ward}) `The canonical height of an algebraic point on an
elliptic curve', \emph{New York J. Math.} \textbf{6} (2000), 331--342.

\bibitem{MR1876275}
(with {\sc M.~Einsiedler} and {\sc T.~Ward})  `Entropy and the canonical height',
  \emph{J. Number Theory} \textbf{91} (2001), no.~2, 256--273.

\bibitem{MR1961589}
(with {\sc T.~Ward})  `Primes in divisibility sequences', \emph{Cubo Mat.
  Educ.} \textbf{3} (2001), no.~2, 245--259.

\bibitem{MR1815962}
(with {\sc M.~Einsiedler} and {\sc T.~Ward})  `Primes in elliptic divisibility
sequences', \emph{LMS J. Comput. Math.} \textbf{4} (2001), 1--13.

\bibitem{MR1938222}
(with {\sc A.~J. van~der Poorten, Y.~Puri}, and {\sc T.~Ward}) `Integer sequences
  and periodic points', \emph{J. Integer Seq.} \textbf{5} (2002), no.~2,
  Article 02.2.3, 10 pp.

\bibitem{7days}
contributed essay in \emph{On the Seventh Day} (John Ashton, ed.),
Master Books (2002).

\bibitem{MR1885618}
(with {\sc I.~Ga{\'a}l, K.~Gy{\"o}ry}, and {\sc C.~R{\"o}ttger})  `On the spatial
  distribution of solutions of decomposable form equations', \emph{Math. Comp.}
  \textbf{71} (2002), no.~238, 633--648.

\bibitem{MR2041076}
(with {\sc P.~Rogers} and {\sc T.~Ward})  `A higher-rank {M}ersenne problem', in
  \emph{Algorithmic number theory ({S}ydney, 2002)}, in \emph{Lecture Notes in
  Comput. Sci.} \textbf{2369}, pp.~95--107 (Springer, Berlin, 2002).

\bibitem{MR2044517}
(with {\sc M.~Einsiedler} and {\sc T.~Ward})  `Morphic heights and periodic points',
  in \emph{Number theory ({N}ew {Y}ork, 2003)}, pp.~167--177 (Springer, New
  York, 2004).

\bibitem{MR1990179}
(with {\sc A.~van~der Poorten, I.~Shparlinski}, and {\sc T.~Ward}) \emph{Recurrence
  sequences}, in \emph{Mathematical Surveys and Monographs} \textbf{104}
  (American Mathematical Society, Providence, RI, 2003).

\bibitem{MR2079832}
(with {\sc M.~Einsiedler} and {\sc T.~Ward})  `Periodic points for good reduction
  maps on curves', \emph{Geom. Dedicata} \textbf{106} (2004), 29--41.

\bibitem{MR2045409}
(with {\sc V.~Miller} and {\sc N.~Stephens}) `Primes generated by elliptic curves',
\emph{Proc. Amer. Math. Soc.} \textbf{132} (2004), no.~4, 955--963.

\bibitem{MR2135478}
(with {\sc T.~Ward}) \emph{An introduction to number theory}, in
  \emph{Graduate Texts in Mathematics} \textbf{232} (Springer-Verlag London
  Ltd., London, 2005).

\bibitem{MR2164113}
(with {\sc H.~King}) `Prime powers in elliptic divisibility sequences',
  \emph{Math. Comp.} \textbf{74} (2005), no.~252, 2061--2071.

\bibitem{MR2200959}
(with {\sc I.~E. Shparlinski})  `Prime divisors of sequences associated to
  elliptic curves', \emph{Glasg. Math. J.} \textbf{47} (2005), no.~1, 115--122.

\bibitem{MR2155503}
(with {\sc K.~Gy{\"o}ry})  `On some arithmetical properties of solutions of
  decomposable form equations', \emph{Math. Proc. Cambridge Philos. Soc.}
  \textbf{139} (2005), no.~1, 27--40.

\bibitem{MR2180241}
(with {\sc V.~Stangoe} and {\sc T.~Ward})  `Orbit counting with an isometric
  direction', in \emph{Algebraic and topological dynamics}, in \emph{Contemp.
  Math.} \textbf{385}, pp.~293--302 (Amer. Math. Soc., Providence, RI, 2005).

\bibitem{MR2220263}
(with {\sc G.~Mclaren} and {\sc T.~Ward}) `Primitive divisors of elliptic
  divisibility sequences', \emph{J. Number Theory} \textbf{118} (2006), no.~1,
  71--89.

\bibitem{MR2309982}
(with {\sc S.~Stevens, D.~Tamsett}, and {\sc T.~Ward})  `Primes generated by
  recurrence sequences', \emph{Amer. Math. Monthly} \textbf{114} (2007), no.~5,
  417--431.

\bibitem{MR2365225}
(with {\sc J.~Reynolds} and {\sc S.~Stevens}) `On the denominators of rational
  points on elliptic curves', \emph{Bull. Lond. Math. Soc.} \textbf{39} (2007),
  no.~5, 762--770.

\bibitem{MR2339472}
(with {\sc R.~Miles, S.~Stevens} and {\sc T.~Ward})  `Orbit-counting in
  non-hyperbolic dynamical systems', \emph{J. Reine Angew. Math.} \textbf{608}
  (2007), 155--182.

\bibitem{MR2429645}
(with {\sc P.~Ingram, V.~Mah{\'e}}, and {\sc S.~Stevens})  `The uniform primality
  conjecture for elliptic curves', \emph{Acta Arith.} \textbf{134} (2008),
  no.~2, 157--181.

\bibitem{MR2428520}
(with {\sc G.~Harman}) `On primitive divisors of {$n^2+b$}', in
  \emph{Number theory and polynomials}, in \emph{London Math. Soc. Lecture Note
  Ser.} \textbf{352}, pp.~142--154 (Cambridge Univ. Press, Cambridge, 2008).

\bibitem{MR2480276}
(with {\sc K.~Eisentr{\"a}ger})  `Descent on elliptic curves and {H}ilbert's
  tenth problem', \emph{Proc. Amer. Math. Soc.} \textbf{137} (2009), no.~6,
  1951--1959.

\bibitem{MR2548983}
(with {\sc V.~Mah{\'e}})  `A generalization of {S}iegel's theorem and
  {H}all's conjecture', \emph{Experiment. Math.} \textbf{18} (2009), no.~1,
  1--9.

\bibitem{MR2486632}
(with {\sc P.~Ingram} and {\sc S.~Stevens})  `Primitive divisors on twists of
  {F}ermat's cubic', \emph{LMS J. Comput. Math.} \textbf{12} (2009), 54--81.

\bibitem{MR2529995}
(with {\sc C.~R{\"o}ttger} and {\sc T.~Ward}) `The continuing story of zeta',
  \emph{Math. Intelligencer} \textbf{31} (2009), no.~3, 13--17.

\bibitem{arXiv}
(with {\sc O.~Phuksuwan} and {\sc S.~Stevens})
`The Uniform Primality Conjecture for the Twisted Fermat Cubic',
{\tt arXiv:1003.2131 [math.NT]} (2010).

\bibitem{MR2550149}
(with {\sc R.~Miles, S.~Stevens} and {\sc T.~Ward})  `Dirichlet series for finite
  combinatorial rank dynamics', \emph{Trans. Amer. Math. Soc.} \textbf{362}
  (2010), no.~1, 199--227.

\bibitem{spectrum}
(with {\sc J.~Griffiths}) `Dual rectangles',
\emph{Math. Spectrum} \textbf{43} (2010/11), 110--114.

\bibitem{MR2915472}
(with {\sc K.~Eisentr{\"a}ger} and {\sc A.~Shlapentokh})
`Hilbert's tenth problem and Mazur's conjectures in complementary subrings of number fields',
\emph{Math. Res. Lett.} \textbf{18} (2011), no.~6, 1141--1162.

\bibitem{MR2826450}
(with {\sc T.~Ward})  `A repulsion motif in {D}iophantine equations',
  \emph{Amer. Math. Monthly} \textbf{118} (2011), no.~7, 584--598.

\end{squarerefs}

%

\affiliationone{Tom Ward\\
The Executive Office\\The Palatine Centre\\Durham University\\DH1 3LE\\
United Kingdom
     \email{t.b.ward@durham.ac.uk}}
\end{document}